\documentstyle[
11pt,amsfonts]{article}

\newcommand {\sph} {{\mathbb S}}
\newcommand {\rp} {{\mathbb R \mathbb P}}
\newcommand {\rel} {{\mathbb R}}
\newcommand {\com} {{\mathbb C}}
\newcommand {\nat} {{\mathbb N}}

\newtheorem{proposition}{Proposition}[section]
\newtheorem{theorem}{Theorem}[section]
\newtheorem{corollaryth}[theorem]{Corollary}

\newtheorem{remark}[proposition]{Remark}

\newtheorem{definition}[proposition]{Definition}

\renewcommand{\theequation}{\mbox{\arabic{section}.\arabic{equation}}}

\newcommand{\nc}{\newcommand}

\newcommand {\pr} {\bf}

\newcommand{\proof} {   \begin{flushright}
                        ///
                        \end{flushright}
                }

\newcommand{\defin} { \hspace*{\fill} $\Box$ }

\newcommand{\foot}[1] { }
\renewcommand{\foot}[1] { \footnote{#1} }

\def	\d	{{ \ {\rm d} }}

\def    \W      {{ {\mathcal W} }}
\nc{\energ}[1]	{{ e_{#1} }}

\def	\poin	{{ poin }}

\def	\geu	{{ g_{euc} }}

\nc{\gpo}[1]	{{ g_{\poin}^{#1} }}
\nc{\gpu}[1]	{{ g_{\poin,{#1}} }}

\nc{\gpmo}[1]	{{ g_{\poin,m}^{#1} }}
\nc{\gpmu}[1]	{{ g_{\poin,m,{#1}} }}
\nc{\tgpmo}[1]	{{ \tilde g_{\poin,m}^{#1} }}

\nc{\gplo}[1]	{{ g_{\poin,\lambda}^{#1} }}

\nc{\conf}[1]	{{ [{#1}] }}

\def    \diff	     {{ \stackrel{\approx}{\longrightarrow} }}

\nc{\doub}[1]{{ \ddot{#1} }}
\nc{\dd}{ \begin{displaymath} }
\nc{\df}{ \end{displaymath} }
\nc{\dcd}{ \begin{displaymath} \begin{array}{c}}
\nc{\dcf}{ \end{array} \end{displaymath} }
\nc{\ee}{ \begin{equation} }
\nc{\ef}{ \end{equation} }
\nc{\ad}{ \begin{array}{c} }
\nc{\af}{ \end{array} }

\begin{document}

\begin{center}
{\huge \bf Rigidity and non-rigidity results for conformal immersions}
\\ \ \\
Tobias Lamm \\
Fakult\"at f\"ur Mathematik\\
Karlsruher Institut f\"ur Technologie (KIT)\\
Englerstra\ss e 2,
D-76131 Karlsruhe, Germany, \\
email: tobias.lamm@kit.edu \\
\ \\
Reiner Michael Sch\"atzle \\
Fachbereich Mathematik der
Eberhard-Karls-Universit\"at T\"ubingen, \\
Auf der Morgenstelle 10,
D-72076 T\"ubingen, Germany, \\
email: schaetz@everest.mathematik.uni-tuebingen.de \\
 
\end{center}
\vspace{1cm}

\begin{quote}

{\bf Abstract:}
In this paper we show a quantitative rigidity result for the minimizer of the Willmore functional among all projective planes in $\rel^n$ with $n\ge 4$. We also construct an explicit counterexample to a corresponding rigidity result in codimension one, by showing that an Enneper surface might split-off during a blow-up process. For conformal immersions of spheres with large enough Willmore energies, we construct explicit counterexamples to a quantitative rigidity result and this complements the recently obtained rigidity results in \cite{nguy.lamm.conf}.
\end{quote}

\vspace{1cm}



\setcounter{equation}{0}

\section{Introduction} \label{intro}

It is a classical result of Codazzi that two-dimensional surfaces in euclidean space which are totally umbilic are parts of a round sphere or a plane. This result was made quantitative in two papers of DeLellis and M\"uller \cite{del.muell}, \cite{del.muell.2}, in which they showed that closed surfaces in $\rel^3$ with small enough tracefree second fundamental form $A^0$ in $L^2$ (or equivalently with Willmore energy close to its absolute minimum value $4\pi$ among all closed surfaces) have to be $W^{2,2}$-close to a round sphere. Additionally, the conformal factor of the pull-back metric has to be $L^\infty$-close to the one of the round sphere. These results relied on delicate estimates for conformal immersions in the Hardy space which were derived by M\"uller and Sverak \cite{muell.sver}. One of the key features of the estimates of DeLellis and M\"uller is that the $W^{2,2}$-estimate for the difference of the immersion and a standard immersion of a round sphere, resp. the $L^\infty$-estimate of the difference of the conformal factors, depends linearly on the $L^2$-norm of $A^0$.

Recently the authors were able to extended these results to surfaces in $\rel^n$, see \cite{lamm.schae.rig}.

Recall that for a smooth immersion $f: \Sigma \to \rel^n\ $ of a closed surface,
we have by the Gau\ss\ equations and the Gau\ss-Bonnet theorem
\begin{displaymath} 
	\W(f)
	= \frac{1}{4} \int \limits_{\Sigma}
	|A_f|^2 \d \mu_f + \pi \chi(\Sigma)
	= \frac{1}{2} \int \limits_{\Sigma}
	|A^0_f|^2 \d \mu_f + 2 \pi \chi(\Sigma),
\end{displaymath}
where
\begin{displaymath}
\W(f)=\frac14 \int \limits_{\Sigma} |H_f|^2 \d \mu_f
\end{displaymath}
is the Willmore energy of $f$. Critical points of $\W$ are called Willmore surfaces.

In the case of spherical surfaces $\Sigma$ and $n=3$, Bryant \cite{bryant84} was able to classify all (smooth) critical points $f_W$ of $\W$. More precisely, he showed that they are inversions of complete minimal surfaces with finite total curvature and embedded planer ends. Additionally he showed that the Willmore energy is quantized in the sense that
\begin{displaymath}
\W(f_W) =4\pi m,
\end{displaymath}
where $m$ is the number of ends of the minimal surface associated to $f_W$ and the values $m=2,3$ are not attained since there are no minimal surfaces with two or three ends satisfying the above conditions. A similar result was shown to be true for $n=4$ by Montiel \cite{montiel}. The result of Bryant was extended to possibly branched Willmore spheres with at most three branch points (including multiplicity) by the first author and Nguyen \cite{nguy.lamm.br}. It was shown that the Willmore energy remains quantized under this assumption. A fact which is no longer true without the restriction on the number of branch points as was observed by Chen and Li \cite{chen.li} and Ndiaye and the second author \cite{ndiaye.schae}. Once singularities are allowed, the energy values $8\pi$ and $12\pi$ show up and they are realized by inversions of the catenoid, resp. the Enneper surface and the trinoid. 

In a recent paper, the first author and Nguyen \cite{nguy.lamm.conf}, were able to obtain quantitative rigidity results for immersions which are close in energy to the inverted catenoid or the inverted Enneper surface and which have at least a multiplicity two point resp. a branch point of branch order two. These immersions have to be $W^{2,2}$-close to the inverted catenoid resp. inverted Enneper surface modulo M\"obius transformations and reparametrizations. Additionally, in higher codimensions, a corresponding quantitative rigidity result was shown to hold for immersions which are close in energy to an inversion of the so called Chen graph (see section 3.1 for more details). These results were obtained by a contradiction argument and hence the linear estimate which was present in the works of DeLellis-M\"uller and the authors of this paper on immersions which are close in energy to the round sphere, is not known so far.

The goal of the present paper is two-fold: First we show that the contradiction argument for the quantitative rigidity result can be extended to non-orientable surfaces and we prove that immersions of $\rp^2$ in $\rel^n$, with $n\ge 4$, whose Willmore energy is close to $6\pi$, the minimal value which is attained by all M\"obius transformations of the stereographic image of the Veronese embedding as was shown by Li and Yau \cite{li.yau}, have to be $W^{2,2}$-close to these surfaces after applying an appropriate M\"obius transformation. Additionally, we show that the conformal factors have to be $L^\infty$-close to each other.

In the second part of the paper we first prove a classification theorem (see Theorem \ref{classificationlimit}) for possible limits modulo M\"obius transformations and reparametrizations of a sequence of possibly branched conformal immersions which converges weakly in $W^{2,2}(\sph^2,\sph^n)$ and in energy to a possibly branched limiting conformal immersion. Our result says that all possible limits have to be either finitely-covered round spheres or a M\"obius transformation and reparametrization of the limit of the original sequence.

Then we use this result in order to show that the above mentioned rigidity results are in some sense optimal. More precisely, we show that conformal immersions of $\rp^2$ in $\rel^3$, which are in energy close to the absolute minimum $\W=12\pi$, do not have to be $W^{2,2}$-close to an immersion attaining the minimum value. The reason for this is a non-compactness property of the moduli space of all such immersions, which was already observed by Bryant \cite{bryant}. Namely, there exists a branched immersion of $\rp^2$ into $\rel^3$ with Willmore energy $12\pi$. The way we construct our counterexample to the rigidity result is by glueing a rescaled version of the Enneper surface into this branched minimizer, thereby desingularizing the immersion. One of the blow-up limits is trivially the Enneper surface itself, which by the above mentioned classification theorem contradicts a possible rigidity result. This answers a question of Bryant \cite{bryant}.

In the case of conformal immersions of spheres in $\rel^n$, $n\ge 4$, we construct counterexamples to possible rigidity results for critical values of the Willmore energy $\W(f)=4\pi m$ once $m\ge 2$. We do this by glueing two appropriately rescaled copies of the conformal minimal immersions $f:\com \to \com^2$, $f(z)=(z^m,z)$ into a $m$-fold branched plane $(z^m,0)$ and then we invert the resulting surface. Together with the classification theorem this yields again a contradiction to a rigidity result. 

In codimension one a similar construction works for all values of $m=2k+1$, $k\in \nat$, by using the so called higher order Enneper surfaces (see section 3.2). For $m=4$ we construct a counterexample by using minimal surfaces classified by Lopez \cite{lopez}.
\newline
In the following we give a brief outline of the paper. 

In section 2 we show the quantitative rigidity result for the stereographic image of the Veronese embedding of $\rp^2$.

In section 3 we present detailed constructions of the above mentioned counterexamples to the quantitative rigidity results.

In the Appendix we show convergence results for conformal factors under appropriate convergence assumptions on the corresponding (branched) conformal immersions. In particular, we establish closeness results for conformal factors which were left open in \cite{nguy.lamm.conf}.

\section{Rigidity for the Veronese embedding}
By a result of Li and Yau (see Theorem 4 in \cite{li.yau}) we know that for all smooth immersions $f:\rp^2 \to \rel^n$ with $n\ge 4$, we have
\begin{displaymath}
\W(f)\ge 6\pi
\end{displaymath}
and equality is attained for all M\"obius transformations of the stereographic image of the Veronese embedding $V:\rp^2 \to \sph^4$, which is given by
\begin{displaymath}
V(x,y,z)=\frac1{\sqrt{3}}(yz,xz,xy,\frac12(x^2-y^2),\frac1{2\sqrt{3}}(x^2+y^2-2z^2)).
\end{displaymath}
Here we want to study a rigidity result related to this fact. Before we can formulate the main result of this section we have to recall the definition of a (branched) $W^{2,2}$-conformal immersion from \cite{kuwert.li}.
\begin{definition}
Let $ \Sigma$ be a Riemann surface. A map $ f \in W^{2,2}_{loc}( \Sigma, \rel^n)$ is called a conformal immersion if in any local conformal coordinates $ (U, z )$, the metric $ g_{i j} = \langle \partial _i f , \partial_j f \rangle $ is given by
\begin{displaymath}
g_{ij} = e^{2u} \delta _{ij}, \quad u \in L^\infty_{loc} (U).
\end{displaymath} 
The set of all $ W^{2,2}$-conformal immersions of $\Sigma$ is denoted $W^{2,2}_{conf} ( \Sigma, \rel^n)$. 

Moreover, a map $ f \in W^{2,2} (\Sigma, \rel^{n})$ is called a branched conformal immersion (with locally square integrable second fundamental form) if $ f \in W ^ {2,2}_ {conf} ( \Sigma \backslash S, \rel^n)$ for some discrete set $ S \subset \Sigma$ and if for each $ p \in S$ there exists a neighbourhood $ \Omega_p $ such that in local conformal coordinates
 \begin{displaymath}
\int_{\Omega_p\backslash \{p\} } |A_f|^2 d\mu_f < \infty.   
\end{displaymath}
Additionally, we either require that $ \mu_{f} (\Omega_p \backslash \{p\} )< \infty $ or that $ p $ is a complete end. 
\end{definition}
\begin{remark}
We note that it follows from Theorem 3.1 in \cite{kuwert.li} that every branch point $p\in \Sigma$ of a $W^{2,2}$-branched conformal immersion $f:\Sigma \to \rel^n$ has a well-defined branch order $m(p)\in \nat_0$.
\end{remark}
Now we are in a position to formulate our rigidity result.
\begin{theorem} \label{intro.theorem-ver}
There exists a number $\delta_0>0$ so that for every $0<\delta < \delta_0$ and every $\ f\in W^{2,2}_{conf}(\rp^2, \rel^n)\ $ with $n\ge 4$ and
\begin{displaymath}
\W(f) \le 6\pi +\delta,
\end{displaymath}
there exists a constant $\omega(\delta)$ with $\omega(\delta)\to 0$ as $\delta \to 0$, a M\"obius transformation $\sigma:\rel^n \to \rel^n$ and a M\"obius transformation of the stereographic image of the Veronese embedding $f_V:\rp^2 \to \rel^4 \subset \rel^n$ so that
\begin{displaymath}
	\| \sigma \circ f - f_V\|_{W^{2,2}(\rp^2,\rel^n)}+ \| u - u_V\|_{L^\infty(\rp^2)} \le \omega(\delta),
\end{displaymath}
where $u$ and $u_V$ are the conformal factors of the conformal immersions $\sigma \circ f$ and $f_V$.
\defin
\end{theorem}

In the proof of this result we will make heavy use of Proposition 4.1 in \cite{kuwert.li} (see also \cite{riviere} for related results). Note that orientability is not an assumption in this result.
\begin{proposition}\label{kuwertli}
Let $\Sigma$ be a closed Riemann surface and $f_k\in W^{2,2}_{conf}(\Sigma, \rel^n)$ be a sequence of conformal immersions satisfying 
\begin{displaymath}
\W(f_k) \le \Lambda < +\infty.
\end{displaymath}
Then for a subsequence there exist M\"obius transformations $\sigma_k$ and a finite set $S\subset \Sigma$ such that
\begin{displaymath}
\sigma_k \circ f_k \to f
\end{displaymath}
weakly in $W^{2,2}_{loc}(\Sigma \backslash S, \rel^n)$, where $f:\Sigma \to \rel^n$ is a branched conformal immersion with square integrable second fundamental form. Moreover, if $\Lambda <8\pi$ then $f$ is unbranched and topologically embedded.
\end{proposition}
\begin{remark}\label{externsionKuwertLi}
We note that the above Proposition extends directly to sequences of branched conformal immersions with a uniform bound on the number of branch points (see e.g. Proposition 3.3 in \cite{nguy.lamm.conf}).
\end{remark}
\begin{remark}\label{bound.conformalfactor}
From the proof of Proposition 4.1 in \cite{kuwert.li} it actually follows that there exists a constant $C<\infty$ so that
\begin{displaymath}
\|u_k\|_{L^\infty_{loc} (\Sigma \backslash S)} \le C,
\end{displaymath}
where $(\sigma_k \circ f_k)^* \delta_{\rel^n} =e^{2u_k} g_0$ and $g_0$ is a smooth background metric on $\Sigma$.

\rm{
This can be seen as follows:

The M\"obius transformations $\sigma_k$ are a composition of a translation (which we can ignore when it comes to bounding the conformal factor), a dilation by a factor $r_k$ so that the conformal factors of $f_k$ have mean value zero on $\Sigma$, and an inversion at a ball of radius $1$ centered at a point $x_0$ whose distance to $r_k^{-1} f_k(\Sigma)$ is bigger than or equal to one.

Now the conformal factors $v_k$ of $r_k^{-1} f_k$ are bounded in $L^\infty_{loc}(\Sigma \backslash S)$ as was shown by Kuwert and Li. Hence one concludes that $\tilde f_k:=r_k^{-1} f_k \rightarrow f$ in
$W^{2,2}_{loc}(\Sigma \backslash  S)$. 

In the case $ \limsup_{k \rightarrow \infty}
\mu_{\tilde f_k}(\Sigma) < \infty$, one concludes $f \in W^{2,2}(\Sigma)$,
and for a subsequence we set $\sigma_k = r_k^{-1}$ and the claim is proved.

For the other alternative $\mu_{\tilde f_k}(\Sigma) \rightarrow \infty$, we choose an inversion as above and get
$I(\tilde f_k) \to I(f)$ weakly in $W^{2,2}_{loc}(\Sigma \backslash S,\rel^n)$.
This can be seen by considering
$U(x) \subset \subset \Sigma \backslash S$ and by noting that $\tilde f_k \rightarrow f$
weakly in $W^{2,2}(U(x),\rel^n)$, in particular the convergence is uniform on $U(x)$.
Hence $|\tilde f_k| \leq R$ is uniformly bounded on $U(x)$ and by the choice of the inversion $I = I_{x_0}$ we have $|f - x_0| \geq 1$ on $\Sigma
\supseteq U(x)$. Since $I$ is smooth on $B_{2R}(0) \backslash B_1(x_0)$, we conclude $I(\tilde f_k) \rightarrow I(f)$ in $W^{2,2}(U(x),\rel^n)$.
More precisely, we have that $|\det DI|$ is bounded from above and below by positive constants on the compact set $\overline{B_{2R}(0)} \backslash B_1(x_0)$ and hence $I^* \delta_{\rel^n} = e^{2v} \delta_{\rel^n}$ where $v \in L^\infty(B_{2R}(0) \backslash B_1(x_0))$. Therefore we have
\begin{displaymath}
\big(I(\tilde f_k)\big)^* \delta_{\rel^n} = (\tilde f_k)^* (e^{2v} \delta_{\rel^n})
        = e^{2 (v_k + v \circ (\tilde f_k))} \delta_{\rel^n},
\end{displaymath}
and the new conformal factor $v_k + v \circ (\tilde f_k)$ remains bounded in $U(x)$, and therefore in $L^\infty (\Sigma \backslash S)$. Here we set $\sigma_k = I(r_k^{-1} \cdot)$ and the claim is again proved.}
\end{remark}
{\pr Proof of Theorem \ref{intro.theorem-ver}:} \\
We argue by contradiction. Hence we assume that there exists a sequence $\delta_k\to 0$, a sequence of $W^{2,2}$-conformal immersions $f_k:\rp^2 \to \rel^n$ satisfying
\begin{displaymath}
6\pi \le \W(f_k) \le 6\pi +\delta_k
\end{displaymath}
and a number $\varepsilon>0$ so that
\begin{displaymath}
	\| \Phi_k \circ f_k - f_V\|_{W^{2,2}(\rp^2,\rel^n)} \ge \varepsilon,
\end{displaymath}
for every sequence of M\"obius transformations $\Phi_k:\rel^n \to \rel^n$ and all M\"obius transformations of the stereographic image of the Veronese embedding $f_V:\rp^2 \to \rel^4\subset \rel^n$.

Using Proposition \ref{kuwertli}, we get the existence of M\"obius
transformations $\sigma_k:\rel^n\to \rel^n$ and an at most finite set $S \subseteq \rp^2$ so that $\tilde f_k:= \sigma_k \circ f_k \rightarrow f$ weakly in
$W^{2,2}_{loc}(\rp^2 \backslash S,\rel^n)$. 
Since 
\begin{displaymath}
\W(f_k) \to 6\pi<8 \pi,
\end{displaymath}
we know additionally that $f\in W^{2,2}_{conf}(\rp^2,\rel^n)$ is a conformal immersion without branch points and hence we have that
\begin{displaymath}
 g := f^* \delta_{\rel^n} = e^{2u} g_{\rp^2}, \ \ \ where \ \ u \in L^\infty(\rp^2).
\end{displaymath}
Using an approximation argument of Schoen-Uhlenbeck as in \cite{kuwert.schaetzle.2013} Proposition 5.2 resp. Theorem 5.1 and using Theorem 4 in \cite{li.yau}, we get that
\begin{displaymath}
\W(f) \geq \inf_{h: \rp^2 \rightarrow \rel^n \,\ smooth} \W(h)
        \geq 6 \pi.
\end{displaymath}
Note again that these results do not require the surface to be orientable.
 
Next, using Gauss-Bonnet and the Gauss equation as in Proposition 5.2 of \cite{kuwert.schaetzle.2013} we conclude
\begin{displaymath}
\int_{\rp^2} |A_f|^2 \d \mu_f = 4\W(f) -4 \pi \chi(\rp^2)
\geq \limsup_{k \rightarrow \infty}
\int_{\rp^2} |A_{\tilde f_k}|^2 \d \mu_{\tilde f_k}=20\pi.
\end{displaymath}
It follows from Remark \ref{bound.conformalfactor} that the conformal factor $\tilde u_k$ of $\tilde f_k$ is bounded in $L^\infty_{loc}(\rp^2 \backslash S)$.

Going back to the proof of Proposition 4.1 in \cite{kuwert.li} and defining
\begin{displaymath}
 \tilde \alpha_k := |A_{\tilde f_k}|^2 \d \mu_{\tilde f_k} \rightarrow \tilde \alpha,
\end{displaymath}
we get
\begin{displaymath}
|A_f|^2 \d \mu_f \leq \tilde \alpha 
\end{displaymath}
in $\rp^2$. We first get this in $\rp^2 \backslash S$, and and then on all of $\rp^2$ since
$\int_S |A_f|^2 \d \mu_f = 0$. 

On the other hand we have that
\begin{displaymath}
\tilde \alpha(\rp^2) \leq \limsup_{k \rightarrow \infty}
\int_{\rp^2} |A_{\tilde f_k}|^2 \d \mu_{\tilde f_k}
\leq \int_{\rp^2} |A_f|^2 \d \mu_f
\end{displaymath}
and therefore
\begin{displaymath}
\tilde \alpha = |A_f|^2 \d \mu_f.
\end{displaymath}
In particular we conclude that
\begin{displaymath}
\tilde \alpha({p}) = 0
\end{displaymath} 
for all $p \in \rp^2$.

Defining a new concentration set 
\begin{displaymath}
\tilde S=\{ p\in \rp^2: \,\ \tilde \alpha (\{ p\}) >0\}
\end{displaymath}
we get $\tilde S = \emptyset$. Note that $\tilde S =S$ if $\sigma_k$ only consists of translations and dilations.

Hence, by repeating the arguments of Proposition 4.1 in \cite{kuwert.li} we get that for the rescaled immersions $\hat f_k = \tilde r_k^{-1} \tilde f_k$, where $\tilde r_k$ is chosen so that the mean value of the conformal factor $\hat u_k=\tilde u_k - \log \tilde r_k$ of $\hat f_k$ is equal to zero, $\hat u_k$ is uniformly bounded in $L^\infty(\rp^2)$.

Since the conformal factors $\tilde u_k$ are already bounded in $L^\infty_{loc}(\rp^2 \backslash S)$ (see Remark \ref{bound.conformalfactor}), 
we conclude that $\log \tilde r_k$ has to be uniformly bounded and therefore $\tilde u_k$ is uniformly bounded in $L^\infty(\rp^2)$.

Arguing as in Proposition 5.3 of \cite{kuwert.schaetzle.2013} or in Remark 2 after Proposition 6.1 in \cite{schae.comp-will12}, we get that $\tilde f_k \rightarrow f$ strongly in $W^{2,2}(\rp^2,\rel^n)$. Using the arguments in \cite{sim.will} one obtains the smoothness of $f$ and by Theorem 4 in \cite{li.yau} we get that $f$ is a M\"obius transformation of the stereographic image of the Veronese embedding, contradicting our assumption. 

The closeness result for the conformal factors now follows from the previous convergence considerations and Theorem \ref{conv:conffac}.
\proof

\section{Non-rigidity results}

In this section we use a glueing argument in order to construct counterexamples to rigidity results for conformal immersions of spheres in arbitrary codimension once the corresponding energy level is high enough.
Our key ingredient in the construction is the following Theorem which classifies possible limits modulo M\"obius transformations and reparametrizations of a strongly converging sequence of conformal immersions.
\begin{theorem}\label{classificationlimit}
Let $f_k\in W^{2,2} (\sph^2,\sph^n)$ be a sequence of possibly branched conformal immersions which converges weakly in $W^{2,2}(\sph^2, \sph^n)$ to a possibly branched conformal immersion $f\in W^{2,2}(\sph^2,\sph^n)$ with singular set $S$. Moreover, we assume that 
\begin{displaymath}
\int_{\sph^2} |A_{f_k}|^2 \d \mu_{f_k} \to \int_{\sph^2} |A_f|^2 \d \mu_f
\end{displaymath}
and that $f\in C^1(\sph^2 \backslash S,\sph^n)$.  

For any sequence of reparametrizations $\varphi_k:\sph^2 \to \sph^2$ and M\"obius transformations $\Phi_k:\sph^n \to \sph^n$ so that $\Phi_k \circ f_k \circ \varphi_k$ converges locally weakly in $W^{2,2}(\sph^2 \backslash S_1,\sph^n)$ (where $S_1$ is again an at most finite set of points) to a possibly branched conformal immersion $h\in W^{2,2}(\sph^2,\sph^n)$, we have that $h$ is either a finitely-covered sphere or $h=\Phi \circ f \circ \varphi$ for a reparametrization $\varphi:\sph^2 \to \sph^2$ and a M\"obius transformation $\Phi:\sph^n \to \sph^n$. 
\end{theorem}
{\pr Proof:} \\
First of all we note that by Fatou's lemma there exists a constant $c_0>0$ so that
\begin{displaymath}
\mu_{f_k}(\sph^2) \ge c_0 \ \ \ and \ \ \  \mu_{\Phi_k \circ f_k}(\sph^2) \ge c_0.
\end{displaymath}
Additionally, we claim that after composing the sequences $f_k$ resp. $\Phi_k \circ f_k$ with an inversion, there exists a constant $C<\infty$ so that
\begin{displaymath}
\mu_{f_k}(\sph^2) , \mu_{\Phi_k \circ f_k}(\sph^2) \le C.
\end{displaymath}
In order to see this, we note that up to a subsequence, $\mu_{f_k} \to \mu$ and $\mu_{\Phi_k \circ f_k} \to \tilde \mu$ weakly as Radon measures. Hence it follows from Proposition A.2 in \cite{lamm.schae.rig} that $f_k(\sph^2) \to spt \mu$ and $\Phi_k(f_k(\sph^2)) \to spt \tilde \mu$ locally in the Hausdorff distance. 
Now there exists $x_0 \in \rel^n \backslash (spt \mu \cup spt \tilde \mu)$, i.e. there exists a $\delta>0$ so that
\begin{displaymath}
d(x_0 , f_k(\sph^2)), \,\ d(x_0, \Phi_k(f_k(\sph^2))) \ge \delta>0.
\end{displaymath}
Inverting both sequences at $B_\delta(x_0)$ we conclude that the images of the two sequences are contained in $B_\delta(x_0)$ and the weak convergence properties remain true. Moreover, the new surfaces satisfy the area bound by \cite{sim.will}.

Next we let $x_k=\Phi_k^{-1}(\infty)$ and we consider three cases:

\underline{Case 1:} $x_k \to \infty$

After a small rotation of $\sph^n$ we can assume that $\Phi_k(\infty)=\infty$, hence
\begin{displaymath}
\Phi_k(x) =\lambda_k O_k x +v_k,
\end{displaymath}
for all $x\in \rel^n$, where $O_k$ is orthogonal. In this case we get the bound
\begin{displaymath}
c_0 / C \le \lambda_k^2 = \mu_{\Phi_k \circ f_k}(\sph^2) / \mu_{f_k} (\sph^2)\le C/ c_0.
\end{displaymath}
Using that $f_k(\sph^2), (\Phi_k \circ f_k)(\sph^2) \subseteq B_\delta(x_0)$, this in turn implies a uniform bound for the $v_k$'s and hence $\Phi_k \to \Phi$ up to a subsequence. In the following we assume without loss of generality that $\Phi=id$. 

By assumption $f_k$ converges to $f$ uniformly on $\sph^2$ and $f_k \circ \varphi_k$ converges to $h$ locally uniformly on $\sph^2 \backslash S_1$. For every $q \in \sph^2 \backslash S_1$ and $\varphi_k(q)\to \tilde q \in \sph^2$ we get after choosing a subsequence
\begin{displaymath}
h(q) \leftarrow (f_k\circ \varphi_k)(q) \rightarrow f(\tilde q)
\end{displaymath}
and by continuity of $h$ this implies $h(\sph^2) \subseteq f(\sph^2)$.

Next, we choose $p_0 \in \sph^2 \backslash (S_1 \cup h^{-1} (f(S)))$ and we note that it follows from the finiteness of $\int_{\sph^2} |A_h|^2 \d \mu_h$ and Corollary 3 in \cite{nguy12} that the set $h^{-1}(f(S))$ is finite. Hence, it follows from Proposition 7.1 in \cite{schae.comp-will12} (or again Corollary 3 in \cite{nguy12}) that $f^{-1}(h(p_0))=\{q^1,\ldots,q^m\}$, $q^i\neq q^j$ for $i\neq j$, is also finite and $q^i\notin S$ for all $i\in \{1,\ldots,m\}$.
Since $f$ is a $C^1$-immersion away from finitely many points in $S$, and by using the inverse function theorem, we additionally get the existence of small balls $B(q^i)$ and $B(h(p_0))$ so that $B(q^i)\cap B(q^j)=\emptyset$ for $i\neq j$, $f=f_i:B(q^i)\diff B(h(p_0))\cap f(B(q^i))$ and
\begin{displaymath}
f^{-1}(B(h(p_0)))= \cup_{i=1}^m B(q^i).
\end{displaymath}
Next we choose another ball $B(p_0)$ so that $h(\overline{B(p_0)})\subseteq B(h(p_0))$ and so that $f_k\circ \varphi_k$ converges uniformly to $h$ on $\overline{B(p_0)}$.  

Now we let $p_k\in \overline{B(p_0)}$ with $p_k\to p \in \overline{B(p_0)}$ and $q_k:= \varphi_k(p_k)\to q\in \sph^2$. We conclude that
\begin{displaymath}
f(q) \leftarrow f_k(q_k)=(f_k\circ \varphi_k)(p_k)\to h(p) \in h(\overline{B(p_0)}) \subseteq B(h(p_0)),
\end{displaymath}
and therefore $q\in f^{-1} (B(h(p_0)))=\cup_{i=1}^m B(q^i)$. In particular this shows that $\varphi_k(p_k)=q_k \in \cup_{i=1}^m B(q^i)$ for $k$ large enough, and hence 
\begin{displaymath}
\varphi_k (\overline{B(p_0)})\subseteq \cup_{i=1}^m B(q^i)
\end{displaymath}
and for $\overline{B(p_0)}$ connected we get the existence of an $i\in \{1,\ldots,m\}$ with
\begin{displaymath}
\varphi_k(\overline{B(p_0)})\subseteq B(q^i).
\end{displaymath}
Using this fact we get $q_k=\varphi_k(p_k)\in B(q^i)$ and hence $q \in \overline{B(q^i)}$, i.e. $q\notin B(q^j)$ for $j\neq i$. On the other hand we know that $q\in \cup_{j=1}^m B(q^j)$ and therefore $q\in B(q^i)$ and $q=f_i^{-1}(h(p))$. In particular this shows that
\begin{displaymath}
\varphi_k |_{\overline{B(p_0)}} \to f_i^{-1} \circ h
\end{displaymath}
uniformly and the right hand side is again in $W^{2,2}\cap W^{1,\infty}(B(p_0))$.

Altogether we get that $\varphi_k \to \varphi$ locally uniformly on $\sph^2 \backslash (S_1\cup h^{-1}(f(S)))$ and $\varphi:\sph^2 \backslash (S_1\cup h^{-1}(f(S))) \to \sph^2$ is in $W^{2,2}_{loc}\cap W^{1,\infty}_{loc}(\sph^2 \backslash (S_1 \cup h^{-1}(f(S))),\sph^2)$ with $h=f\circ \varphi$.

Since both $f$ and $h$ are conformal on $\sph^2 \backslash (S_1 \cup h^{-1}(f(S)))$, it follows that $\varphi$ is also conformal and we can choose an orientation on the connected set $\sph^2 \backslash (S_1 \cup h^{-1}(f(S)))$ so that $\varphi$ is holomorphic on this set. 

For $\sph^2 \backslash (S_1\cup h^{-1}(f(S))) \ni p_k \to p \in S_1\cup h^{-1}(f(S))$ and $\varphi(p_k) \to q$ we get
\begin{displaymath}
f(q) \leftarrow (f\circ \varphi) (p_k) = h(p_k) \rightarrow h(p)
\end{displaymath}
and thus $q\in f^{-1}(h(S_1 \cup h^{-1}(f(S))))\neq \sph^2$. Therefore the points in $S_1 \cup h^{-1}(f(S))$ are removable singularities for $\varphi$ by the Theorem of Casorati-Weierstrass and we obtain a holomorphic extension $\varphi:\sph^2 \to \sph^2$ with $h=f \circ \varphi$.
Since $\varphi$ is holomorphic and non-constant it is a branched immersion and we get
\begin{displaymath}
\int_{\sph^2} |A_f|^2 \d \mu_f =\lim_{k\to \infty} \int_{\sph^2} |A_{f_k}|^2 \d \mu_{f_k} \ge \int_{\sph^2} |A_h|^2 \d \mu_h = \deg \varphi \int_{\sph^2} |A_f|^2 \d \mu_f
\end{displaymath}
and hence we conclude that $\deg \varphi =1$, which implies via the Riemann-Hurwitz formula (see e.g. Theorem I.2.7 in \cite{far.kra} or Theorem 2.5.2 in \cite{jo.rie}) that $\varphi$ is a diffeomorphism.

Therefore we conclude that $h$ is a reparametrization of $f$ and $\int_{\sph^2} |A_h|^2 \d\mu_h=\int_{\sph^2} |A_f|^2 \d \mu_f $.

\underline{Case 2:} $x_k \to x\in \rel^n$ and $x\not\in f(\sph^2)$

We let $I_{x_k}(x):= (x-x_k) /|x-x_k|^2$ and define $\Psi_k:= \Phi_k \circ I_{x_k}^{-1}$. Then $\Psi_k (\infty)=\infty$ and again we have $\Psi_k =\lambda_k O_k +v_k$. Now $I_{x_k}$ converges smoothly to $I_x$ in a neighborhood of $f(\sph^2)$, in particular $I_{x_k} \circ f_k \to I_x \circ f$ and we can replace $f_k$ by $I_{x_k} \circ f_k$, $\Phi_k$ by $\Psi_k$ and we are back in case 1).

\underline{Case 3:} $x_k \to x\in \rel^n$ and $x \in f(\sph^2)$

We let $I_{x_k}$ and $\Psi_k$ be as in case 2), and we note that
\begin{displaymath}
\mu_{I_{x_k} \circ f_k}(\sph^2) \to \infty
\end{displaymath}
and hence
\begin{displaymath}
\lambda_k^2 = \mu_{\Psi_k \circ I_{x_k} \circ  f_k}(\sph^2) / \mu_{I_{x_k} \circ f_k} (\sph^2)\to 0.
\end{displaymath}
After choosing a subsequence we can also assume that $v_k \to v \in \rel^n \cup \{\infty\}$. For $y_k \to y \in \sph^n \backslash \{x\}$ we get 
\begin{displaymath}
I_{x_k}(y_k) \to I_x(y) \in \rel^n
\end{displaymath}
and 
\begin{displaymath}
\Phi_k (y_k) =\Psi_k ( I_{x_k} (y_k)) =\lambda_k O_k I_{x_k} (y_k) +v_k \to v.
\end{displaymath}
Choosing a metric on $\rel^n \cup \{\infty\} \cong \sph^n$, we conclude for every $\delta>0$
\begin{displaymath}
\Phi_k^{-1} (\sph^n \backslash B_{\delta} (v)) \subseteq B_\delta (x)
\end{displaymath}
for all $k$ large enough.

Using the conformal invariance of $\int |A^0|^2 \d \mu$, this implies that for every $\delta>0$ 
\begin{displaymath}
\int_{h^{-1}(\rel^n \backslash B_{2\delta} (v)) } |A_h^0|^2 \d \mu_h \le \liminf_{k\to \infty} \int_{ (\Phi_k \circ f_k)^{-1}(\rel^n \backslash B_\delta(v))} |A_{\Phi_k \circ f_k}^0|^2 \d \mu_{\Phi_k \circ f_k}
\end{displaymath}
\begin{displaymath}
\le \liminf_{k\to \infty} \int _{f_k^{-1}(B_\delta(x))} |A_{f_k}^0|^2 \d \mu_{f_k} .
\end{displaymath}
The last term converges to zero as $\delta \to 0$, since by the assumptions of the theorem no local energy concentration is possible for the sequence of immersions $f_k$. It follows that $A_h^0\equiv 0$ and $h$ is a finitely-covered sphere.
\proof
\begin{remark}\label{doublecover}
Using the oriented double cover $\pi:\sph^2 \to \rp^2$ it is easy to see that the above Lemma extends directly to the case of sequences of possibly branched conformal immersions $f_k\in W^{2,2}(\rp^2,\sph^n)$.
\end{remark}

\subsection{Higher codimensions}
The situation of conformal immersions $f:\sph^2 \to \rel^n$ with $n\ge 4$ is easier to handle, since there exists a counterexample at a lower energy level than in the codimension one case.

In this setting we have a conformal minimal immersion $f_C:\com \to \com^2$ with one end of multiplicity two, the so called Chen graph. It can be parametrized by
\begin{displaymath}
f_C(z)= (z^2,z)
\end{displaymath}
and it satisfies $\int_{\com} |A_{f_C}|^2 \d \mu_C =4\pi$. Inverting the Chen graph at a point $x_0\notin f_C(\com)$, yields a $W^{2,2}$-branched conformal immersion $f_c:\sph^2 \to \rel^4$ with exactly one branch point $p\in \sph^2$ of branch order $m(p)=1$. In particular, it follows from the generalized Gauss-Bonnet theorem (see Corollary 2 and 3 in \cite{nguy12}) that
\begin{displaymath}
\W(f_c)=8\pi \ \ \ and \ \ \ \int_{\sph^2} |A_{f_c}|^2  \d \mu_c = 20\pi.
\end{displaymath}
Moreover, any branched conformal immersion $f_b:\sph^2 \to \rel^n$, $n\ge 4$, with at least one branch point of branch order $m(p)=1$ satisfies
\begin{displaymath}
\W(f_b)\ge 8\pi \ \ \ and \ \ \ \int_{\sph^2} |A_{f_b}|^2  \d \mu_b \ge 20\pi.
\end{displaymath}
The following result was shown in Theorem 4.4 of \cite{nguy.lamm.conf}
\begin{theorem}\label{rigidity4}
There exists $\delta_0>0$ such that for every $0<\delta<\delta_0$ and every immersion $f\in W^{2,2}_{conf,br}(\sph^2, \rel^n)$, $ n \geq 4$, where $ f $ has exactly one branch point $p\in \sph^2$ of branch order $m(p) = 1 $ and which satisfies $20 \pi \leq  \int_{\sph^2} |A|^2 d\mu \le 20\pi+\delta$, there exists a M\"obius transformation $\sigma:\rel^n\to \rel^n$, a reparametrization $\phi: \sph^2 \to  \sph^2$, a
constant $\omega(\delta)$, with $\omega(\delta)\to 0$ as $\delta \to 0$, and a standard immersion $f_{c}\in W^{2,2}_{conf,br}(\sph^2,\rel^n)$ of an inverted Chen graph with
\begin{displaymath}
||\sigma \circ f \circ \phi-f_{c}||_{W^{2,2}( \sph^2,\rel^n)} \le \omega(\delta). \label{chen_close}
\end{displaymath}
\end{theorem}
Our goal here is to show that there is no rigidity result available for conformal immersions of the sphere into $\rel^n$, $n\ge 4$, with at least one double point. This is in sharp contrast to Theorem 1.3 in \cite{nguy.lamm.conf}, where a rigidity result for immersions from the sphere with exactly one double point has been proved in the case $n=3$.

We note that every $f\in W^{2,2}_{conf} (\sph^2,\rel^n)$ with at least one double point satisfies (see Corollary 2 and 3 in \cite{nguy12})
\begin{displaymath}
\W(f)\ge 8\pi \ \ \ and \ \ \ \int_{\sph^2} |A_{f}|^2  \d \mu_f \ge 24\pi
\end{displaymath}
and the immersions attaining equality have been characterized in \cite{hof.oss}.

In the following we construct a sequence of smooth conformal immersions $f_\varrho:\sph^2 \to \rel^4$, $\varrho \to 0$, which have at least one double point and which satisfy $\W(f_\varrho)\searrow 8\pi$, $\int_{\sph^2} |A_{f_\varrho}|^2 \d \mu_{f_\varrho} \searrow 24 \pi$ but so that for any sequence of M\"obius transformations $\Phi_\varrho:\rel^4 \to \rel^4$ and any sequence of reparametrizations $\psi_\varrho:\sph^2 \to \sph^2$, the new sequence $\Phi_\varrho \circ f_\varrho \circ \psi_\varrho$ cannot converge weakly in $W^{2,2}(\sph^2,\rel^4)$ to a limiting conformal immersion $f:\sph^2 \to \rel^4$ with exactly one double point and with $\W(f)=8\pi$ and $\int_{\sph^2} |A_f|^2 \d\mu_f=24\pi$.

Here is the precise construction: 
\newline
For every $\varrho>0$ we define an immersion $f_\varrho: \com \to \com^2$ by
\begin{displaymath}
f_\varrho(z):= \varrho^2 f_C(\varrho^{-1} z)= (z^2,\varrho z).
\end{displaymath}

Next we choose a smooth function $\varphi:\com \to \com$ with $\varphi(z)=z$ for all $|z|\le \frac12$, $\varphi \in C^\infty_c(B_1(0))$ and we define a new immersion $h_\varrho:\com \to \com^2$ by
\begin{displaymath}
h_\varrho (z)=(z^2, \varrho \varphi(z) ).
\end{displaymath}
Note that by this construction we are glueing a rescaled Chen graph into a double plane and moreover we have
\begin{displaymath}
\W (h_\varrho) \to 0 \ \ \ and \ \ \ \int_{\com} |A_{h_\varrho}|^2 \d \mu_{h_\varrho}\to 4\pi
\end{displaymath}
as $\varrho \to 0$. The immersion $h_\varrho$ has a complete end of multiplicity two.

Now we define the inversion at the unit sphere $I:\rel^4 \to \rel^4$, $I(x)=\frac{x}{|x|^2}$ and we let
\begin{displaymath}
\tilde h_\varrho(w) := I ( h_\varrho( 1/\bar w) ).
\end{displaymath}
The so defined immersion $\tilde h_\varrho: \com 
\rightarrow \rel^4$ is a smooth immersion away from $0$ which satisfies (see again Corollary 2 and 3 in \cite{nguy12} and note that $h_\varrho(1/ \bar w)$ has a complete end of multiplicity two at the origin and a branch point of branch order one at infinity with $h_\varrho(0)=0$)
\begin{displaymath}
\W (\tilde h_\varrho) \to 4\pi \ \ \ and \ \ \ \int_{\com} |A_{\tilde h_\varrho}|^2 \d \mu_{\tilde h_\varrho}\to 12 \pi
\end{displaymath}
as $\varrho \to 0$.

By construction we have $\tilde h_\varrho(w) = (w^2,0)$ for $0<|w| \leq 1$, and hence we can repeat the above steps in order to glue in another copy of the rescaled Chen graph $f_\varrho$ in $B_1\subset \com$. In this way we obtain an immersion $\bar h_\varrho:\com \to \com^2$ satisfying 
\begin{displaymath}
\W (\bar h_\varrho) \to 4\pi \ \ \ and \ \ \ \int_{\com } |A_{\bar h_\varrho}|^2 \d \mu_{\bar h_\varrho}\to 16\pi
\end{displaymath}
as $\varrho \to 0$ and every point in $\bar h_\varrho (\partial B_1(0))$ is a double point.

We note that for every $z\in B_{(2\varrho)^{-1}}(0)$ 
\begin{displaymath}
\varrho^{-2}  \bar h_\varrho(\varrho z)=(z^2, z)
\end{displaymath}
and hence a subsequence of $\varrho^{-2} \bar h_\varrho (\varrho \cdot)$ converges locally weakly in $W^{2,2}$ on all of $\com$ to a Chen graph as $\varrho \to 0$.

Finally, in order to get a compact image surface, we note that we can find a ball $B_1(x_0)$ with $x_0\in \rel^4$ so that $\varrho^{-2} \bar h_\varrho (\com) \cap B_1(x_0)=\emptyset$ for all $\varrho>0$ and hence, by defining $\hat h_\varrho :\com \cup \{\infty\} \cong \sph^2 \to \rel^4$, $\hat h_\varrho (z)=I_{x_0} (\varrho^{-2} \bar h_\varrho (\varrho z))$, we obtain a sequence of immersions satisfying
\begin{displaymath}
\W (\hat h_\varrho) \to 8\pi \ \ \ and \ \ \ \int_{\sph^2} |A_{\hat h_\varrho}|^2 \d \mu_{\hat h_\varrho}\to 24\pi.
\end{displaymath}
and which converges weakly in $W^{2,2}_{loc}(\sph^2 \backslash \{N\},\rel^4)$ to an inverted Chen graph.

If we now assume that there exists a sequence of M\"obius transformations $\Phi_\varrho:\rel^4 \to \rel^4$ and a sequence of reparametrizations $\varphi_\varrho:\sph^2 \to \sph^2$ so that $\Phi_\varrho \circ \hat h_\varrho \circ \varphi_\varrho$ has a subsequence which converges weakly in $W^{2,2}(\sph^2,\rel^4)$ to a conformal immersion $h\in W^{2,2}(\sph^2,\rel^4)$ with at least one double point and so that
\begin{displaymath}
\W ( h) = 8\pi \ \ \ and \ \ \ \int_{\sph^2} |A_{h}|^2 \d \mu_{h}= 24\pi,
\end{displaymath}
then it follows that $h$ has exactly one double point $d\in \sph^2$ and $h\in C^\infty_{loc} (\sph^2 \backslash h^{-1}(d),\rel^4)$ and this yields a contradiction to Theorem \ref{classificationlimit}. 

In order to see that $h$ has exactly one double point, we note that the existence of a second double point would imply that we can invert the image surface $h(\sph^2)\subset \rel^4$ at one of the double points in order to obtain a minimal surface $h_1$ with two ends of multiplicity one, $\int |A_{h_1}|^2 \d  \mu_{h_1} =8\pi$ and which has at least one other double point, without loss of generality we assume it is at the origin. In this case we can invert $h_1$ at the unit ball in $\rel^4$ and we obtain again a minimal surface $\tilde h_1$ (see e.g. Corollary 2 in \cite{nguy12}). Using the transformation formula for the mean curvature under inversions we get
\begin{displaymath}
0= H_{\tilde h_1} = |h_1|^4 (H_{h_1} + 4 h_1 ^\perp |h_1|^{-2}) =4 |h_1|^{2} h_1^\perp,
\end{displaymath}
where $(v)^\perp$ denotes the normal component of $v\in \rel^4$ and hence we conclude that $h_1 \subseteq T h_1$. In particular we conclude that $h_1$ is a cone and therefore $A_{h_1} \equiv 0$, which contradicts the fact that $\int | A_{h_1}|^2 \d \mu_{h_1} =8\pi$. 

We conclude that $h_1$ is a conformal minimal immersion and hence it is smooth. This shows that $h$ is smooth in $\sph^2$ away from the preimages of its double point.

The same construction can be performed with $f_C$ being replaced by the conformal minimal immersion $f_m(z)=(z^m,z)$ for every $m\ge 2$. Note that $\int_{\com} K_{f_m} \d \mu_{f_m}=2\pi(1-m)$ and therefore $\int_{\com} |A_{f_m}|^2 \d \mu_{f_m} = 4\pi (m-1)$. Glueing two appropriately rescaled immersions $f_m$ into a $m$-plane $(z^m,0)$ and inverting the resulting surface, yields a sequence of immersions $f_\varrho:\sph^2 \to \com^2$ with $\W(f_\varrho)\to 4m \pi$ and $\int_{\sph^2} |A_{f_\varrho}|^2 \d \mu_{f_\varrho} \to 8\pi(2m-1)$, as $\varrho \to 0$. Note that these energy values are also attained by inverting minimal surfaces with $m$-ends of multiplicity one. Hence, we obtain a non-rigidity result for all these energy levels.

In particular we obtain the following Theorem.
\begin{theorem}\label{nonrigidityhigher}
For any number $m\in \nat$, $m\ge 2$, there exists a sequence of smooth conformal immersions $f_\varrho:\sph^2 \to \rel^4$, $\varrho \to 0$, which have at least one point of multiplicity $m$ and which satisfy $\W(f_\varrho)\searrow 4m \pi$, $\int_{\sph^2} |A_{f_\varrho}|^2 \d \mu_{f_\varrho} \searrow (2m-1) 8 \pi$ but so that for any sequence of M\"obius transformations $\Phi_\varrho:\rel^4 \to \rel^4$ and any sequence of reparametrizations $\psi_\varrho:\sph^2 \to \sph^2$, the new sequence $\Phi_\varrho \circ f_\varrho \circ \psi_\varrho$ cannot converge weakly in $W^{2,2}(\sph^2,\rel^4)$ to a limiting conformal immersion $f:\sph^2 \to \rel^4$ with exactly one point of multiplicity $m$ and with $\W(f)=4m \pi$ and $\int_{\sph^2} |A_f|^2 \d\mu_f=(2m-1)8 \pi$.
\end{theorem}

\subsection{Codimension one}
Next we look at conformal immersions $f:\sph^2 \to \rel^3$ and we note that in this case we have a surface with one end of multiplicity three, the Enneper surface
\begin{displaymath}
f_E(z)= -\frac19 (z^3,0) +\frac13 (x, -y,x^2-y^2),
\end{displaymath}
where $z=x+iy$. This is a conformal minimal immersion which satisfies $\int_{\com} |A_{f_E}|^2 \d \mu_{f_e} =8\pi$. Inverting the Enneper surface yields a $W^{2,2}$-branched conformal immersion $f_e:\sph^2 \to \rel^3$ with exactly one branch point $p\in \sph^2$ of branch order $m(p)=2$. In particular, it follows from Corollary 2 and 3 in \cite{nguy12} that
\begin{displaymath}
\W(f_e)=12\pi \ \ \ and \ \ \ \int_{\sph^2} |A_{f_e}|^2  \d \mu_{f_e} = 32\pi.
\end{displaymath}

Next we note that every $f\in W^{2,2}_{conf} (\sph^2,\rel^3)$ with at least one triple point satisfies (see again Corollary 2 and 3 in \cite{nguy12})
\begin{displaymath}
\W(f)\ge 12\pi \ \ \ and \ \ \ \int_{\sph^2} |A_{f}|^2  \d \mu_f \ge 40\pi
\end{displaymath}
and the inversion of a trinoid is an example of a conformal immersion with a triple point attaining equality in both estimates. Note that after an inversion at a triple point of an immersion attaining equality in the above estimates, one obtains a complete minimal immersion with $\int_\com |A|^2 \d \mu =16\pi$ and these have been classified by Lopez \cite{lopez}.

As in the previous subsection we construct a sequence of conformal immersions with $W(f_k)\to 12\pi$ and $\int_{\sph^2} |A_{f_\varrho}|^2 \d \mu_{f_\varrho} \to 40\pi$ but for which no composition with a M\"obius transformation resp. reparametrization converges weakly to an immersion $f:\sph^2 \to \rel^3$ with at least one triple point and whose Willmore energy is equal to $12\pi$ and with $\int_{\sph^2} |A_{f}|^2  \d \mu_f = 40\pi$. 

In order to do this, we define for every $\varrho >0$ the rescaled Enneper surfaces ($z=x+iy$)
\begin{displaymath}
f_\varrho (z)= \varrho^3 f_E(\varrho^{-1}z)= -\frac19 (z^3,0) +\frac13 (\varrho^2 x, -\varrho^2 y, \varrho(x^2-y^2)).
\end{displaymath}
Next we chose a smooth function $\varphi:\com \to \rel$ with $\varphi(z) =1$ for all $|z|\le \frac12$, $\varphi \in C^\infty_c(B_1(0))$ and we define new immersions
\begin{displaymath}
h_\varrho (z)= -\frac19 (z^3,0) + \frac13 \varphi(z) (\varrho^2 x, -\varrho^2 y, \varrho(x^2-y^2)).
\end{displaymath}
As $\varrho \to 0$ we have that 
\begin{displaymath}
\W(h_\varrho) \to 0 \ \ \ and \ \ \ \int_{\com} |A_{h_\varrho}|^2 \d \mu_{h_\varrho} \to 8\pi.
\end{displaymath}
As in the previous subsection we invert this immersion at the unit sphere in both the domain and the image and we glue in another copy of the rescaled Enneper surface in the unit ball. In this way we obtain an immersion $\bar h_\varrho :\com \to \rel^3$ with at least one triple point,
\begin{displaymath}
\W(\bar h_\varrho) \to 8\pi \ \ \ and \ \ \ \int_{\com} |A_{\bar h_\varrho}|^2 \d \mu_{\bar h_\varrho} \to 32\pi
\end{displaymath}
and so that $\varrho^{-3}\bar h_\varrho (\varrho \cdot)$ converges locally weakly in $W^{2,2}$ on $\com$ to an Enneper surface as $\varrho \to 0$. 

After an inversion $I_{x_0}$ at a unit ball around a point $x_0 \notin \bar h_\varrho$ we then get the existence of a sequence of immersions $\tilde h_\varrho:\sph^2 \cong \com \cup \{ \infty \} \to \rel^3$, $\tilde h_\varrho (z)=I_{x_0}(\varrho^{-3} \bar h_\varrho (\varrho z))$ with at least one triple point,
\begin{displaymath}
\W(\tilde h_\varrho) \to 12 \pi \ \ \ and \ \ \ \int_{\sph^2} |A_{\tilde h_\varrho}|^2 \d \mu_{\tilde h_\varrho} \to 40\pi
\end{displaymath}
and which converges weakly in $W^{2,2}_{loc}(\sph^2 \backslash \{N\},\rel^3)$ to an inverted Enneper surface.

Using the same arguments as in the case of higher codimensions we combine this construction with Theorem \ref{classificationlimit} in order to get a contradiction to a rigidity result for conformal immersions $f:\sph^2 \to \rel^3$ with at least one triple point.

This construction can again be extended to a non-rigidity result for inversions of minimal surfaces with $(2m+1)$-ends, $m\in \nat$, of multiplicity one. Namely, there exist the so called higher order Enneper surfaces $f_{HE}:\com \to \rel^3$,
\[
f_{HE}(z)=Re \left( z-(2m+1)^{-1} z^{2m+1}, i(z+(2m+1)^{-1} z^{2m+1}), 2(m+1)^{-1} z^{m+1} \right)
\]
which are conformal minimal immersions with one end of multiplicity $2m+1$. By glueing two rescaled versions of these higher order Enneper surfaces into a $(2m+1)$-plane, we get again a non-rigidity result for all conformal immersions $f:\sph^2 \to \rel^3$ with at least one point of multiplicity $2m+1$ and $\W(f)=(8m+4)\pi$, $\int_{\sph^2} |A_f|^2 \d \mu_f = (32m+8)\pi$.
Note that this construction only works for odd multiplicities, since there do not exist complete minimal surfaces in $\rel^3$ with finite total curvature and one end of even multiplicity by a result of Osserman \cite{osser}. 

An example of a conformal immersion with a quadruple point and minimal Willmore energy resp. $\int |A|^2$ is the Morin surface $f_{Mo}:\sph^2 \to \rel^3$. It satisfies $\W(f_{Mo})=16\pi$ and $\int_{\sph^2} |A_{f_{Mo}}|^2 \d \mu_{f_{Mo}}=56\pi$. In this case we can construct a counterexample to the rigidity as follows: It was shown by Lopez \cite{lopez} that there exists a complete minimal surface $f_1:\com \backslash \{0\}\to \rel$ with one end of multiplicity three and one end of multiplicity one, i.e. with $\int_{\com} K_{f_1} \d \mu_{f_1}=-8\pi$ or $\int_{\com} |A_{f_1}|^2 \d \mu_{f_1}=16\pi$. The surface is given explicitely by its Weierstrass representation and one easily calculates that the blow-down of the surface converges to the union of a triple plane and a single plane,  which intersect in a line, away from the origin. Hence, inverting $f_1$ at a point not in its image (say $0\notin f_1(\com \backslash \{0\}$), we get a branched conformal immersion $f_1^\star:\sph^2 \to \rel^3$ with $f_1^\star(\pm N)=0$. Close to $N$ we can now glue in a triple plane and close to $-N$ we glue in a single plane without changing the energies too much. Additionally, we can assume that close to $N$ the immersion $f_1^\star$ parametrizes an annulus $A_1 := B_1(0) \backslash \overline{B_{1/2}(0)}$ in the triple plane $P_1$ and close to $-N$, $f_1^\star$ parametrizes an annulus $A_2 := B_1(0) \backslash \overline{B_{1/2}(0)}$ in the single plane $P_2$. The set $P_1\cap P_2$ contains a line $L$ in $\rel^3$ and hence the new immersion has quadruple points on $L\cap B_1 \backslash \overline{B_{1/2}(0)}$. 
Finally, we can glue in a rescaled copy of the Enneper surface into the triple plane we already glued in close to $N$. In this way we obtain a sequence of immersions $f_\varrho:\sph^2 \to \rel^3$ with at least one quadruple point and $\W(f_\varrho) \to 16\pi$, $\int_{\sph^2} |A_{f_\varrho}|^2 \d \mu_{f_\varrho} \to 56\pi$. In particular, this construction implies again a non-rigidity result for the Morin surface.

Combining the constructions we obtain the following theorem.
\begin{theorem}\label{nonrigiditythree}
For any number $m\in \nat$, $m=4$ or $m=2k+1$ with $k\in \nat$, there exists a sequence of smooth conformal immersions $f_\varrho:\sph^2 \to \rel^3$, $\varrho \to 0$, which have at least one point of multiplicity $m$ and which satisfy $\W(f_\varrho)\searrow 4m \pi$, $\int_{\sph^2} |A_{f_\varrho}|^2 \d \mu_{f_\varrho} \searrow (2m-1) 8 \pi$ but so that for any sequence of M\"obius transformations $\Phi_\varrho:\rel^3 \to \rel^3$ and any sequence of reparametrizations $\psi_\varrho:\sph^2 \to \sph^2$, the new sequence $\Phi_\varrho \circ f_\varrho \circ \psi_\varrho$ cannot converge weakly in $W^{2,2}(\sph^2,\rel^3)$ to a limiting conformal immersion $f:\sph^2 \to \rel^3$ with exactly one point of multiplicity $m$ and with $\W(f)=4m \pi$ and $\int_{\sph^2} |A_f|^2 \d\mu_f=(2m-1)8 \pi$.
\end{theorem}

\subsection{$\rp^2$ in codimension one}

In this subsection we show that a rigidity result for minimizers $f:\rp^2\to \rel^3$ of $\W$ similar to Theorem \ref{intro.theorem-ver} fails in codimension one. 

It follows from a topological result of Banchoff \cite{banchoff} that every smooth immersion $f:\rp^2\to \rel^3$ has to have at least one triple point, and hence it follows from Theorem 10, Corollary 2 and 3 in \cite{nguy12} that for every conformal immersion $f\in W^{2,2}(\rp^2,\rel^3)$ one has
\begin{displaymath}
\W(f)\ge 12\pi \ \ \ and \ \ \ \int_{\rp^2} |A_f|^2 \d \mu_f \ge 44\pi.
\end{displaymath}
It was shown by Bryant \cite{bryant} and Kusner \cite{kusner2,kusner} that there exists a minimizing conformal immersion, the Boy's surface, $f_B:\rp^2\to \rel^3$ with
\begin{displaymath}
\W(f_B)= 12\pi \ \ \ and \ \ \ \int_{\rp^2} |A_{f_B}|^2 \d \mu_{f_B} = 44\pi.
\end{displaymath}
We note that this surface is an inversion of a complete non-oriented, non-embedded minimal immersion with three ends of multiplicity one. 

On the other hand, Meeks \cite{meeks} constructed a minimal immersion $f_M:\rp^2 \backslash \{p\} \to \rel^3$ with one end $p\in \rp^2$ of mulitiplicity three and $\int K_{f_M} \d \mu_{f_M} =-6\pi$.
The Weierstrass data of this non-orientable minimal surface is (see \cite{lopezmartin}) 
\[
(M,I,g,\eta)= (\com \backslash \{0\}, I(z)=-1/\bar z, g(z)=z^2\frac{z+1}{z-1},\eta= i\frac{(z-1)^2}{z^4} \d z).
\]
Integrating this data shows that for $|z|>>1$ the conformal minimal immersion $\tilde f_M:\com \backslash \{0\}\to \rel^3$ satisfies
\[
\tilde f_M(z)= -(i z^3/6, 0) +( O(z^2),O(z)).
\]
In particular the same expansion remains true for $f_M$ close to the single end $p$. This shows that by performing a blow-down $f_\varrho(z) =  \varrho^3 f_M(\varrho^{-1} z)$ the immersions converge locally smoothly away from the point $p$ to the triple plane $-(i z^3/6,0)$.

Inverting the surface $f_M$ at $\partial B_1(0)$, without loss of generality we assume $0 \notin f_M(\rp^2 \backslash \{p\})$, one obtains a branched conformal immersion $f_m:\rp^2 \to \rel^3$ with a single branch point $p\in \rp^2$ of branch order $m(p)=2$ and which satisfies
\begin{displaymath}
\W(f_m)= 12\pi \ \ \ and \ \ \ \int_{\rp^2} |A_{f_m}|^2 \d \mu_{f_m} = 36\pi.
\end{displaymath}
In the following we glue a rescaled Enneper surface into a rescaling of the branched minimizer $f_m$ in order to construct a sequence of immersions of $\rp^2$ into $\rel^3$ which contradicts a rigidity result similar to Theorem \ref{intro.theorem-ver}. Moreover, the construction shows additionally that minimizing sequences of $W^{2,2}(\rp^2,\rel^3)$-conformal immersions for the Willmore functional may degenerate in the sense that after composing the sequence with suitably chosen sequences of M\"obius transformations and reparametrizations the new sequence can converge to an Enneper surface. This gives an answer to a question mentioned on page 240 of \cite{bryant}.

In order to perform this construction we note that it follows from the above discussion that in a small neigborhoud, say $B_{2\delta}(p) \backslash B_\delta(p)$, around $p$, the immersion $f_m$ looks like a triple plane $-(i z^3/6,0)$ up to an error which smoothly converges to zero as $\delta \to 0$. Hence we can use a similar construction as in subsection $3.2$ and we can glue in a rescaled version $f_\varrho$ of the Enneper surface in $B_\delta(p)$. 

With this construction one obtains a sequence of immersions $\hat f_\varrho :\rp^2 \to \rel^3$ so that 
\begin{displaymath}
\W(\hat f_\varrho) \to 12\pi, \ \ \ \ \int_{\rp^2} |A_{\hat f_\varrho}|^2 \d \mu_{\hat f_\varrho} \to 44\pi
\end{displaymath}
and whose blow-up (in local conformal coordinates around $p$) $\varrho^3 \hat f_\varrho (\varrho^{-1}\cdot): B_{2\delta \varrho}(0)\to \rel^3$ converges locally weakly in $W^{2,2}(\com,\rel^3)$ to the Enneper surface.

Arguing as in the previous two subsections one can combine this construction with Theorem \ref{classificationlimit} (see also Remark \ref{doublecover}) to get the non-rigidity result for the Boy's surface. 
\begin{theorem}\label{nonrigidityenneper}
There exists a sequence of smooth conformal immersions $f_\varrho:\rp^2 \to \rel^3$, $\varrho \to 0$, which have at least one point of multiplicity $3$ and which satisfy $\W(f_\varrho)\searrow 12 \pi$, $\int_{\rp^2} |A_{f_\varrho}|^2 \d \mu_{f_\varrho} \searrow  44\pi$ but so that for any sequence of M\"obius transformations $\Phi_\varrho:\rel^3 \to \rel^3$ and any sequence of reparametrizations $\psi_\varrho:\rp^2 \to \rp^2$, the new sequence $\Phi_\varrho \circ f_\varrho \circ \psi_\varrho$ cannot converge weakly in $W^{2,2}(\rp^2,\rel^3)$ to a limiting conformal immersion $f:\rp^2 \to \rel^3$ with exactly one point of multiplicity $3$ and with $\W(f)=12 \pi$ and $\int_{\rp^2} |A_f|^2 \d\mu_f= 44 \pi$.
\end{theorem}
{\LARGE \bf Appendix}

\begin{appendix}
\renewcommand{\theequation}{\mbox{\Alph{section}.\arabic{equation}}}


\setcounter{equation}{0}
\section{Convergence of conformal factors under strong $W^{2,2}$-convergence of conformal immersions} 
In this appendix we prove a convergence result for the conformal factors of a strongly converging sequence of conformal immersions.
\begin{theorem}\label{conv:conffac}
Let $\Sigma$ be a closed surface with a smooth metric $g_\Sigma$ and let $f_k\in W^{2,2}_{conf}(\Sigma, \rel^n)$ be a sequence of conformal immersions with induced metrics $g_k=f_k^\star \geu= e^{2u_k} g_\Sigma$. We assume that 
\[
f_k\to f\in W_{conf}^{2,2}(\Sigma, \rel^n) 
\]
converges strongly in $W^{2,2}(\Sigma,\rel^n)$ and that there exists a constant $C<\infty$ so that 
\[
\| u_k\|_{L^\infty(\Sigma)} \le C.
\]
Then we have that
\[
\|u_k - u\|_{L^\infty(\Sigma)} \to 0,
\]
where $u$ is the conformal factor of $f$, i.e. $g=f^\star \geu= e^{2u} g_\Sigma$.
\end{theorem}
{\pr Proof:} \\
The assumptions of the theorem imply that $g_k \rightarrow g$
weakly in $W^{1,2}$ and strongly in $L^p$ for all $1\le p <\infty$. In particular, after selecting a subsequence, we conclude that $g_k \rightarrow g$ and $u_k \to u$ pointwise almost everywhere.

We let $\varepsilon >0$ and by using the strong convergece of the $f_k$ in $W^{2,2}$ we find finitely many conformal charts $\varphi_i: U_i
\subseteq \Sigma \diff B_2(0)$, such that the subsets $V_i = \varphi_i^{-1}(B_{1}(0))$ form a covering of $\Sigma$, i.e. $\Sigma = \cup_i V_i$,
and so that for all $i$ we have
\begin{displaymath}
\sup \limits_{k}\int \limits_{U_i} |A_{f_k}|^2 \d \mu_{g_k}
< \varepsilon.
\end{displaymath}
We fix $i$ and define the conformal immersions
\[
\tilde f_{k,i}
:= f_{k} \circ \varphi_i^{-1}: B_2(0) \rightarrow \rel^n
\]
and 
\[
\tilde f_i:= f\circ \varphi_i^{-1}: B_2(0) \rightarrow \rel^n.
\]
Moreover, we let $\tilde g_{k,i} := \tilde f_{k,i}^* \geu = e^{2 \tilde u_{k,i}} \geu$ resp. $\tilde g_{i} := \tilde f_{i}^* \geu = e^{2 \tilde u_{i}} \geu$ and we observe that
\[
(u_{k} - u) \circ \varphi_i^{-1}
= \tilde u_{k,i} - \tilde u_i,
\] 
and in particular this yields $\tilde u_{k,i} \rightarrow
\tilde u_i$ almost everywhere in $B_2(0)$.

Using \cite{schae.comp-will12}, Proposition 5.1, or \cite{kuwert.li}, Corollary 2.4, there exist functions
$v_{k,i}: B_2(0) \rightarrow \rel$ and $v_i:B_2(0) \rightarrow \rel$, so that $\tilde u_{k,i} - v_{k,i}$ and $\tilde u_i -v_i$ are harmonic in $B_2(0)$ and so that
\[
\parallel v_{k,i} \parallel_{L^\infty(B_2(0))}
\leq C_n \int \limits_{B_2(0)} |A_{\tilde f_{k,i}}|^2 \d \mu_{\tilde g_{k,i}}
\leq C_n \varepsilon
\]
resp.
\[
\parallel v_i \parallel_{L^\infty(B_2(0))}\le C_n \int \limits_{B_2(0)} |A_{\tilde f_i}|^2 \d \mu_{\tilde g_i}
\leq C_n \varepsilon.
\]
Using these facts we conclude
\begin{displaymath}
\parallel u_k - u \parallel_{L^\infty(V_i)}
= \parallel \tilde u_{k,i} - \tilde u_i \parallel_{L^\infty(B_{1}(0))}
\leq \| \tilde u_{k,i} - v_{k,i} - (\tilde u_i - v_i) \|_{L^\infty(B_{1}(0))}+C_n \varepsilon
\end{displaymath}
\begin{displaymath}
\leq C \| \tilde u_{k,i} - v_{k,i} - (\tilde u_i - v_i) \|_{L^1(B_2(0))}
+ C_n \varepsilon\\
\leq C \| \tilde u_{k,i} - \tilde u_i \|_{L^1(B_2(0))} + C_n \varepsilon
\rightarrow C_n \varepsilon,
\end{displaymath}
and therefore
\begin{displaymath}
\limsup \limits_{k \rightarrow \infty}
\parallel u_k - u \parallel_{L^\infty(\Sigma)} \leq C_n \varepsilon
\end{displaymath}
from which we get that $u_k \rightarrow u$ in $L^\infty(\Sigma)$.
\proof
As a Corollary of this result we obtain
\begin{corollaryth}\label{conv:conffacbranch}
Let $\Sigma$ be a closed surface with a smooth metric $g_\Sigma$ and let $f_k\in W^{2,2}_{conf,loc}(\Sigma \backslash \{p\}, \rel^n)$ be a sequence of branched conformal immersions with at most one branch point $p\in \Sigma$ of branch order $m$, for some $m \in \nat_0$, and with induced metrics $g_k=f_k^\star \geu= e^{2u_k} g_\Sigma$. We assume that $f_k\to f$ converges weakly in $W^{2,2}(\Sigma,\rel^n)$ and that $f$ is a branched conformal immersion which has again a branch point $p\in \Sigma$ of branch order $m$. Moreover, we assume that 
\[
\int_\Sigma |A_{f_k}|^2 \d \mu_{f_k} \to \int_\Sigma |A_f|^2 \d \mu_f.
\]
Then we have that
\[
\|u_k - u\|_{L^\infty(\Sigma)} \to 0,
\]
where $u$ is the conformal factor of $f$, i.e. $g=f^\star \geu= e^{2u} g_\Sigma$.
\end{corollaryth}
{\pr Proof:} \\
As in the proof of the Theorem, the assumptions of the theorem imply that $g_k \rightarrow g$
weakly in $W^{1,2}$ and strongly in $L^p$ for all $1\le p <\infty$. In particular, after selecting a subsequence, we conclude that $g_k \rightarrow g$ and $u_k \to u$ pointwise almost everywhere.

We let $\varepsilon >0$ and by using the convergence assumptions on $f_k$ and the $L^2$-norm of the second fundamental form of $f_k$, we find finitely many conformal charts $\varphi_i: U_i
\subseteq \Sigma \diff B_2(0)$, such that the subsets $V_i = \varphi_i^{-1}(B_{1}(0))$ form a covering of $\Sigma$, i.e. $\Sigma = \cup_i V_i$,
and so that for all $i$ we have
\begin{displaymath}
\sup \limits_{k}\int \limits_{U_i} |A_{f_k}|^2 \d \mu_{f_k}
< \varepsilon.
\end{displaymath}
Without loss of generality we assume that $p\in V_1$ and $p \notin V_i$, $i\ge 2$.

It follows from Remark \ref{bound.conformalfactor} and the proof of Theorem 3.2 in \cite{nguy.lamm.conf} that there exists a constant $C<+\infty$ so that 
\[
\|u\|_{L^\infty(V_i)}\le C \ \ \ for \ \ i\ge 2
\]
and 
\[
\|u-m\log |\cdot| \|_{L^\infty(V_1)}\le C.
\]

We fix $i$ and define as before the conformal immersions
\[
\tilde f_{k,i}
:= f_{k} \circ \varphi_i^{-1}: B_2(0) \rightarrow \rel^n
\]
and 
\[
\tilde f_i:= f\circ \varphi_i^{-1}: B_2(0) \rightarrow \rel^n.
\]
Moreover, we let $\tilde g_{k,i} := \tilde f_{k,i}^* \geu = e^{2 \tilde u_{k,i}} \geu$ resp. $\tilde g_{i} := \tilde f_{i}^* \geu = e^{2 \tilde u_{i}} \geu$ and we observe that
\[
(u_{k} - u) \circ \varphi_i^{-1}
= \tilde u_{k,i} - \tilde u_i.
\] 
In particular this yields $\tilde u_{k,i} \rightarrow
\tilde u_i$ almost everywhere in $B_2(0)$ for $i\ge 2$ and $\tilde u_{k,i}-m\log |\cdot| \rightarrow
\tilde u_i-m\log |\cdot|$ almost everywhere in $B_2(0)$ for $i=1$.

For $i\ge 2$ we argue as in the proof of Theorem \ref{conv:conffac} in order to get
\begin{displaymath}
\parallel u_k - u \parallel_{L^\infty(V_i)}
= \parallel \tilde u_{k,i} - \tilde u_i \parallel_{L^\infty(B_{1}(0))}
\leq \| \tilde u_{k,i} - v_{k,i} - (\tilde u_i - v_i) \|_{L^\infty(B_{1}(0))}+C_n \varepsilon
\end{displaymath}
\begin{displaymath}
\leq C \| \tilde u_{k,i} - v_{k,i} - (\tilde u_i - v_i) \|_{L^1(B_2(0))}
+ C_n \varepsilon\\
\leq C \| \tilde u_{k,i} - \tilde u_i \|_{L^1(B_2(0))} + C_n \varepsilon
\rightarrow C_n \varepsilon.
\end{displaymath}
In the case $i=1$ we argue as in the proof of Theorem 3.1 in \cite{kuwert.li} in order to construct functions $v_{k,1}: B_2(0) \rightarrow \rel$ and $v_1:B_2(0) \rightarrow \rel$, so that $\tilde u_{k,1} -m\log | \cdot |- v_{k,1}$ and $\tilde u_1-m\log |\cdot| -v_1$ are harmonic in $B_2(0)$ and so that
\[
\parallel v_{k,1} \parallel_{L^\infty(B_2(0))}
\leq C_n \int \limits_{B_2(0)} |A_{\tilde f_{k,1}}|^2 \d \mu_{\tilde f_{k,1}}
\leq C_n \varepsilon
\]
resp.
\[
\parallel v_1 \parallel_{L^\infty(B_2(0))}\le C_n \int \limits_{B_2(0)} |A_{\tilde f_1}|^2 \d \mu_{\tilde f_1}
\leq C_n \varepsilon.
\]
The same argument as in the Theorem above then also yields
\begin{displaymath}
\parallel u_k - u \parallel_{L^\infty(V_1)}\le \parallel (\tilde u_{k,1}-m\log |\cdot|) - (\tilde u_1 -m\log |\cdot|) \parallel_{L^\infty(B_{1}(0))}\le C_n \varepsilon.
\end{displaymath}
Combining the estimates for all $V_i$ implies
\begin{displaymath}
\limsup \limits_{k \rightarrow \infty}
\parallel u_k - u \parallel_{L^\infty(\Sigma)} \leq C_n \varepsilon
\end{displaymath}
from which we get that $u_k \rightarrow u$ in $L^\infty(\Sigma)$.
\proof
\begin{remark}
We note that Theorem \ref{conv:conffac} and Corollary \ref{conv:conffacbranch} show that the rigidity results obtained in Theorem 1.3, 1.4 , 4.2 and 4.4 in \cite{nguy.lamm.conf} directly imply a corresponding rigidity result for the conformal factors in the $L^\infty$-norm as well.
\end{remark}
\end{appendix}


\end{document}